\documentclass[12pt]{article}

\usepackage{amsmath, amssymb, amsthm}
\usepackage{hyperref}
\usepackage{fullpage}
\usepackage{graphicx}

\newcommand{\Z}{{\mathbb Z}}
\newcommand{\F}{{\mathbb F}}
\newcommand{\T}{\rule{0pt}{2.4ex}}
\newcommand{\B}{\rule[-1.2ex]{0pt}{0pt}}

\title{Scrambling Connections Puzzles with Finite Fields}

\author{Alon Danai, Joshua Kou, Andy Latto, Haran Mouli, and Jim Propp}

\date{March 10, 2026}

\begin{document}

\maketitle

\begin{abstract}
\noindent
We define a {\it grid graph} $G$ as a Cartesian product of path-graphs 
$P_n$ or cycle-graphs $C_n$ as shown in Figure~\ref{fig:grids},
and we ask, when can the edge set of a complete graph be expressed as 
a disjoint union of graphs isomorphic to $G$? That is, we are asking
for which grid graphs a $G$-design exists, where a $G$-design is defined 
as a decomposition of a complete graph into edge-disjoint subgraphs 
isomorphic to $G$. We show that when $n$ is an odd prime or the square 
of an odd prime, the toroidal grid-graph $G = C_n \square C_n$ admits 
a $G$-design. In the less symmetrical case of products of path-graphs,
we prove that $G = P_3 \square P_3$ does not admit a $G$-design but 
that $G = P_4 \square P_4$ does. This last result is the special case 
that motivated the present paper: a $P_4 \square P_4$-design corresponds 
to a way of successively scrambling a Connections puzzle so that each 
pair of words occurs adjacently exactly once. Our constructions use 
the arithmetic of finite fields.
\end{abstract}

% \keywords?
% \keywords{graph decompositions, graph designs,
% Cartesian products, finite fields}

\vspace{10px}



\begin{figure}
\begin{center}
\includegraphics[width=4in]{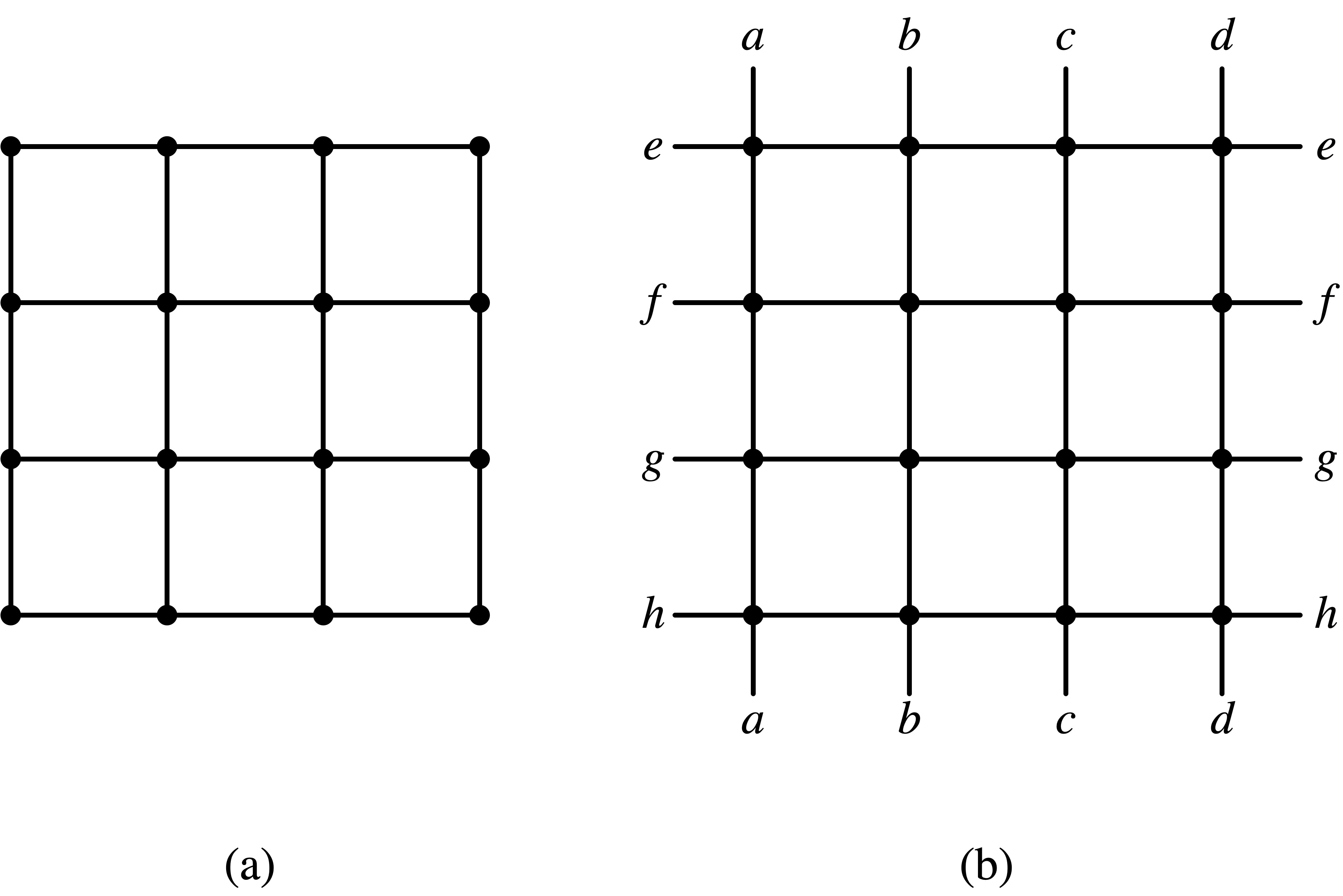}
\end{center}
\caption{Panel (a) shows a plane grid of order 4, isomorphic to $P_4 \square P_4$.
Panel (b) shows a torus grid of order 4, isomorphic to $C_4 \square C_4$;
edges in the torus grid marked with the same label are identified with each other.}
\label{fig:grids}
\end{figure}

\section{Introduction} \label{sec:intro}

One of the oldest problems in graph theory, the {\it probl\`eme de ronde},
predates the general concept of graphs by half a century. In volume 2 of 
his series of books on recreational mathematics~\cite{L}, published around 
1890, Edouard Lucas posed the problem of finding a schedule for seating $2n+1$ 
diners at a round table on each of $n$ successive nights in such a manner 
that each diner sits next to each other diner exactly once. 

\begin{figure}
\begin{center}
\includegraphics[width=4in]{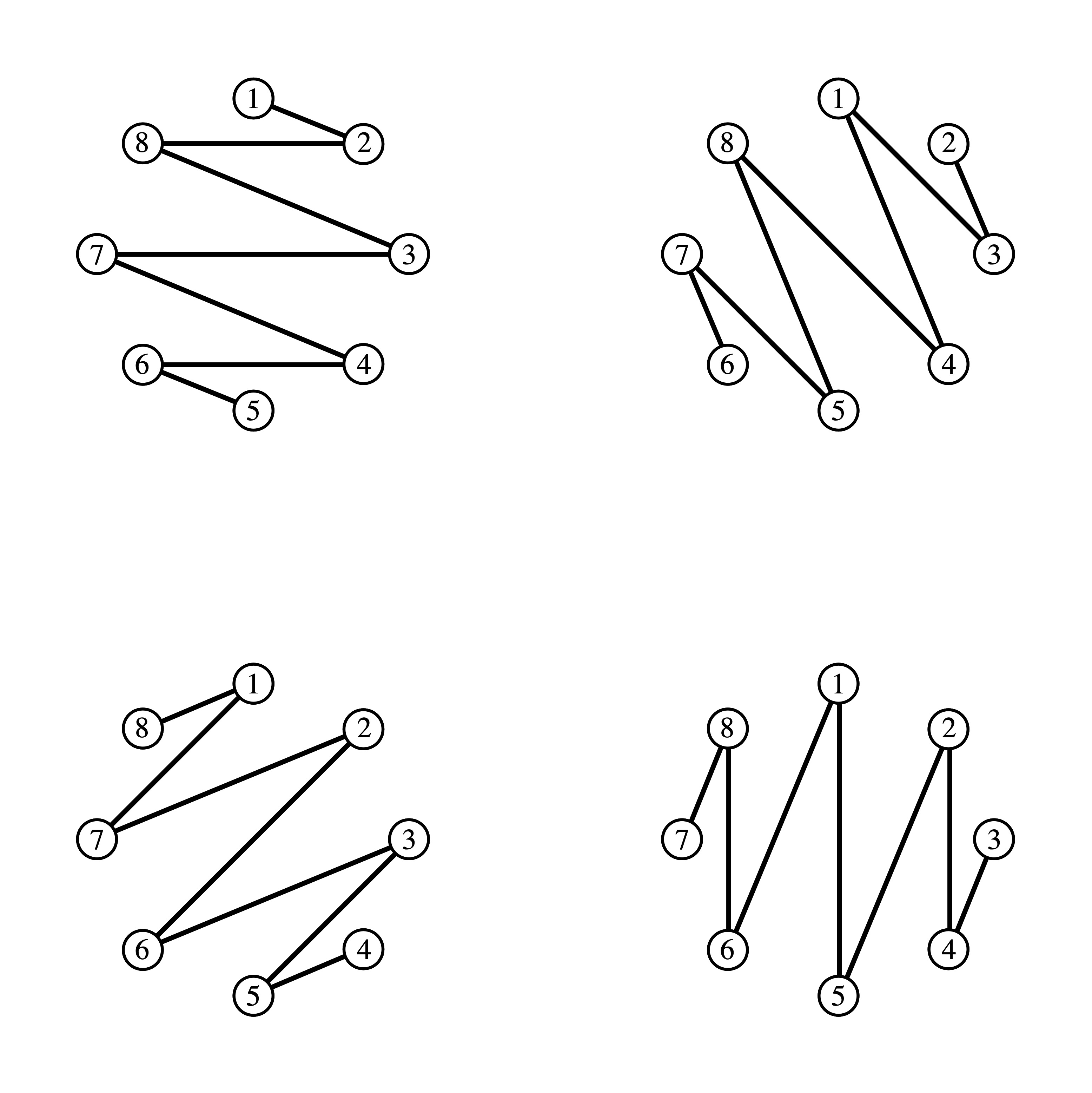}
\end{center}
\caption{A decomposition of $K_8$ into four $P_8$ graphs.}
\label{fig:even}
\end{figure}

Lucas gave an elegant solution that he attributed to F\'elix Walecki.
Walecki's method, slightly modified, can be applied to a related
problem one might call the {\it probl\`eme de banc}, in which we must
find a schedule for seating $2n$ people along a bench on each of $n$ 
successive days in such a manner that each person sits next to each 
other person exactly once. In fact, the solution to the second problem 
is slightly simpler, and is illustrated in the case $n=4$ in 
Figure~\ref{fig:even}. The Figure shows four zig-zag paths which,
when superimposed, give the complete graph on 8 vertices, with no
repeated edges. Rotating each zigzag by 45 degrees clockwise
gives its clockwise neighbor in the Figure. Each zigzag provides 
an ordering to use along the bench.  In graph theory terms, we are 
decomposing the edge-set of $K_8$ (the complete graph on 8 vertices) 
into four copies of $P_8$ (the path graph on 8 vertices).

\begin{figure}
\begin{center}
\includegraphics[width=4in]{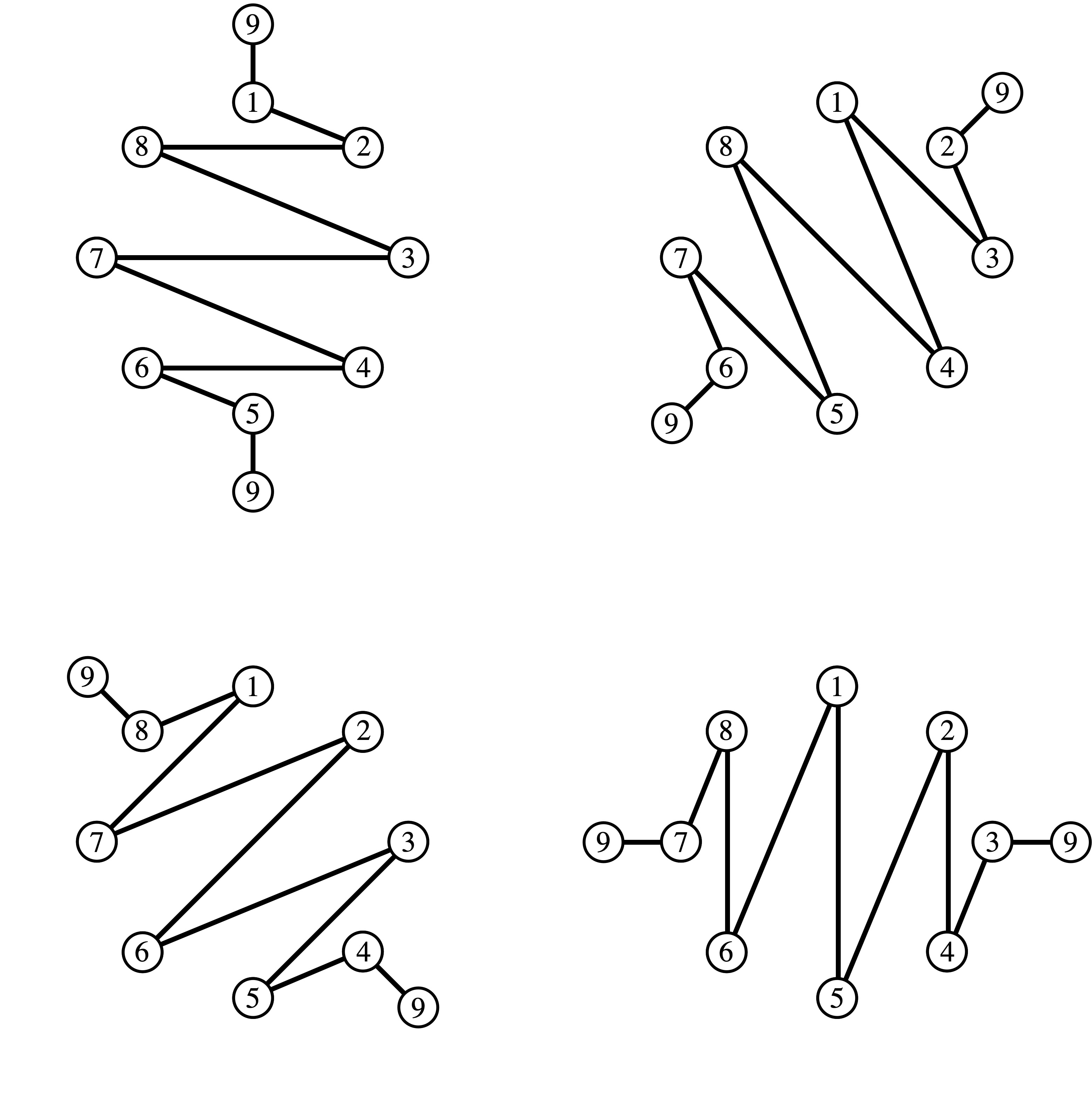}
\end{center}
\caption{A decomposition of $K_9$ into four $C_9$ graphs.
Vertices marked with the same label are identified with each other.}
\label{fig:odd}
\end{figure}

Walecki's solution to the actual {\it probl\`eme de ronde}
is shown, in graph-theoretic form, in Figure~\ref{fig:odd},
for $n=4$.  The $2n+1$st diner adjoined to a solution to the 
{\it probl\`eme de banc}, sitting next to a different pair
of diners each time. In graph theory terms, we are 
decomposing the edge-set of $K_9$ 
into four copies of $C_9$ (the cycle graph on 9 vertices).

Here we extend Lucas' question by looking for $G$-designs
where $G$ is not a cycle graph or a path graph 
but a {\em Cartesian product} of such graphs, where 
$G \square G'$, the Cartesian product of the graphs 
$G=(V,E)$ and $G'=(V',E')$, is defined as the graph
whose vertices are of the form $(v,v')$
with $v \in V$ and $v' \in V'$ 
and whose edges are of the form $\{(v,v'), (w,w')\}$
with either (1) $\{v,w\} \in E$ and $v'=w'$
or (2) $\{v',w'\} \in E$ and $v=w$.
For instance, the graph shown in Figure~\ref{fig:grids}(a)
is $P_4 \square P_4$
while the graph shown in Figure~\ref{fig:grids}(b)
is $C_4 \square C_4$.
We are not aware of any prior work on this extension.

The motivation comes from the clue-scrambling feature
of the electronic version of Connections puzzles,
as introduced by the New York Times in 2023.
Figure~\ref{fig:puzzle} shows a sample Connections puzzle,
too easy to be typical of the genre but illustrative of the core idea.
\begin{figure}
$$
\renewcommand{\arraystretch}{2.5} % increase row height, centered
\begin{array}{|c|c|c|c|}
\hline
\text{OVER} & \text{\!MOUTH\!} & \text{\:HAND\:} & \text{NOSE} \\[1ex]
\hline
\text{HORSE} & \text{SPITE} & \text{MILK} & \text{BITE} \\[1ex]
\hline
\text{FEEDS} & \text{CRY} & \text{GIFT} & \text{CUT} \\[1ex]
\hline
\text{LOOK} & \text{FACE} & \text{YOU} & \text{\:SPILT\:} \\[1ex]
\hline
\end{array}
$$
\caption{A sample Connections puzzle.}
\label{fig:puzzle}
\end{figure}
Solving the puzzle involves noticing that the sixteen words 
can be partitioned into four groups, each of which comes from 
a different English proverb indicating things you shouldn't do:
BITE the HAND that FEEDS YOU,
LOOK a GIFT HORSE in the MOUTH,
CUT off your NOSE to SPITE your FACE,
or CRY OVER SPILT MILK.
In this instance, the four foursomes are related by 
a meta-theme (``things you shouldn't do'');
in actual Connections puzzles, that usually isn't the case,
since such extra structure makes a puzzle too easy solve.

Note that the grid contains the similar words 
``bite'' and ``cut'' next to each other, 
suggesting membership in a foursome of actions 
that divide an object into pieces.
Likewise the proximity of ``mouth'', ``nose'', and ``hand''
suggests shared membership in a foursome of body parts.
These are false leads, and such red herrings are 
commonly included on purpose to make the puzzle trickier.
In the on-line version of the puzzle, there's a button
that allows you to scramble the words
to make the false leads less salient
and to bring out the true patterns,
or at least to make it possible for your brain
to find new connections between words.
The computer code that scrambles the words
uses a pseudorandom number generator
to scramble the sixteen words uniformly;
for instance, there is a 1 out of 16!\ chance
that pressing the Scramble button will have no effect at all.

Empirically one finds that it takes between 20 and 30 scramblings
for all the possible adjacencies between the 16 words to appear.
That is, if one scrambles only 20 times, it is nearly certain that
some pair of words will not have occurred adjacently,
whereas if one scrambles 30 times, it is nearly certain that
every pair of words will have occurred adjacently at least once.

Since each 4-by-4 grid (the original grid and the successive scrambled versions)
has 24 adjacent pairs, and since $16 \choose 2$ is exactly five times 24,
one might hope that if one is extremely lucky,
it might happen that after a mere four scramblings,
each of the possible pairs of words will have occurred adjacently exactly once.
The question is, can this actually happen?
We will show that the answer is {\sc Yes}; that is, 
$K_{16}$ can be partitioned into five copies of $P_4 \square P_4$.
Our proof makes use of the finite field with 16 elements,
though the first construction of such a partition 
was not constructed algebraically but 
was found by Ed Kirkby using computer search~\cite{MO}.
Just as Walecki's way of decomposing $K_{2n+1}$ into copies of $C_n$
has a cyclic symmetry of order $n$,
Kirkby's way of decomposing $K_{16}$ into copies of $P_4 \square P_4$
has a cyclic symmetry of order $5$.

Likewise, since a 3-by-3 grid graph has 12 adjacent pairs of cells,
and since $9 \choose 2$ is exactly three times 12,
one might hope that if the words in a 3-by-3 grid 
are scrambled and then scrambled again in just the right way,
each of the possible pairs of words will occur adjacently exactly once.
Can this actually happen?
A combinatorial argument shows that in this case the answer is {\sc No}; 
that is, $K_{9}$ cannot be partitioned into three copies of $P_3 \square P_3$.

More generally, whenever $n$ is congruent to 0 or 3 mod 4,
so that $n^2 \choose 2$ is a multiple of $2n(n-1)$,
one may ask, for each such $n$, whether $K_{n^2}$
can be partitioned into copies of $P_n \square P_n$.
We have not made much of a dent in this problem,
but we have made some progress with
a more symmetrical version of the problem
in which $P_n \square P_n$ is replaced by $C_n \square C_n$.
Since $n^2 \choose 2$ is a multiple of $2n^2$
whenever $n$ is congruent to 1 or 3 mod 4,
one may ask, for each such $n$, whether $K_{n^2}$
can be partitioned into copies of $C_n \square C_n$;
we will show that this is indeed the case whenever $n$ 
is an odd prime raised to the first or second power.

This article is organized as follows.
In Section~\ref{sec:torus},
we use finite field methods to show that 
when $p$ is an odd prime, $K_{p^2}$ can be partitioned 
into $(p^2-1)/4$ copies of $C_p \square C_p$ (Theorem 1), 
and that $K_{p^4}$ can be partitioned 
into $(p^4-1)/4$ copies of $C_{p^2} \square C_{p^2}$ (Theorem 2).
In Section~\ref{sec:ordinary}, we show that
$K_{9}$ cannot be partitioned into 3 copies of $P_3 \square P_3$ (Theorem 3)
and by modifying the method of Section~\ref{sec:torus} we show that
$K_{16}$ can be partitioned into 5 copies of $P_4 \square P_4$ (Theorem 4).
Section~\ref{sec:comments} presents comments and open questions.

\section{Torus Grids} \label{sec:torus}

We can picture the torus graph $C_n \square C_n$ either 
by modding out the infinite 4-regular graph $\Z \square \Z$
by the translations $(n,0)$ and $(0,n)$
or by drawing the ordinary grid $P_n \square P_n$
and adding wraparound edges along the boundary.
Since $C_n \square C_n$ has $2n^2$ edges
while $K_{n^2}$ has $n^2 (n^2-1)/2$ edges,
it is clear that if $K_{n^2}$ can be decomposed into
disjoint copies of $C_n \square C_n$,
the number of copies must be $(n^2-1)/4$,
so $n$ must be odd.

\bigskip

{\bf Theorem 1:} When $p$ is an odd prime,
$K_{p^2}$ can be partitioned into $(p^2-1)/4$ copies of $C_p \square C_p$.

{\bf Proof:} Associate the vertices of $K_{p^2}$ 
with the elements of $\F_{p^2}$,
which we model as $\F_p[x]/(f(x))$
where $f$ is some primitive quadratic polynomial over $\F_p$,
so that $x$ is a generator of the multiplicative group $\F_{p^2}^*$,
so that $x^{(p^2-1)/2}=-1$.  For all $\alpha, \beta \in \F_{p^2}$
that form a basis for $\F_{p^2}$ over $\F_p$,
let $G(\alpha, \beta)$ be the subgraph of $K_{p^2}$
in which two vertices are joined by an edge
if and only if the corresponding field elements
differ by $\pm \alpha$ or $\pm \beta$.
It is easy to see that each such subgraph
is isomorphic to $C_p \square C_p$. 
The two-element sets $\{ x^0, -x^0 \}$, $\{ x^1, -x^1 \}$, \dots
$\{ x^{(p^2-3)/2}, -x^{(p^2-3)/2} \}$
form a partition of $\F_{p^2}^*$.
Then the graphs $G(x^0,x^1)$, $G(x^2,x^3)$,
\dots, $G(x^{(p^2-5)/2},x^{(p^2-3)/2})$
form the desired partition of $K_{p^2}$.
\qed

\bigskip

We move on to the $p^2$-by-$p^2$ torus.

\bigskip

{\bf Theorem 2:} When $p$ is an odd prime,
$K_{p^4}$ can be partitioned into 
$(p^4-1)/4$ copies of $C_{p^2} \square C_{p^2}$.

{\bf Proof:} 
Associate the vertices of $K_{p^4}$ 
with the elements of $\F_{p^4}$,
and let $x$ be a generator 
of the multiplicative group $\F_{p^4}^*$,
so that $x^{(p^4-1)/2}=-1$.
For all $\alpha, \beta,\gamma,\delta \in \F_{p^4}$
that form a basis for $\F_{p^4}$ over $\F_p$,
let $G(\alpha, \beta,\gamma,\delta)$ be the subgraph of $K_{p^4}$
in which two vertices are joined by an edge
if and only if the corresponding field elements
differ by $\pm \alpha,\pm \beta,\pm\gamma,\pm\delta$.
It is easy to see that each such subgraph
is isomorphic to $C_p \square C_p \square C_p\square C_p$. 
By Kotzig~\cite{Kotzig}, $C_p \square C_p$ can be decomposed into 
two Hamiltonian cycles $C_{p^2}$ partitioning the edges. 
Hence, performing the decomposition on the grouped pairs in
$(C_p \square C_p) \square (C_p\square C_p)$
we obtain four graphs 
$G_1(\alpha, \beta,\gamma,\delta), G_2(\dots), G_3(\dots), G_4(\dots)$, 
each isomorphic to $C_{p^2} \square C_{p^2}$, which 
together decompose $G(\alpha, \beta,\gamma,\delta)$.
The two-element sets $\{ x^0, -x^0 \}$, $\{ x^1, -x^1 \}$, 
\dots $\{ x^{(p^4-3)/2}, -x^{(p^4-3)/2} \}$
form a partition of $\F_{p^4}^*$.
Then the graphs 
$$G(x^0,x^1,x^2,x^3),\ G(x^4,x^5,x^6,x^7),\ \dots, \ G(x^{(p^4-9)/2},x^{(p^4-7)/2},x^{(p^4-5)/2},x^{(p^4-3)/2})$$
form a decomposition of $K_{p^4}$, 
and replacing each $G$ with $G_1,G_2,G_3,$ and $G_4$ 
gives a decomposition in terms of $C_{p^2} \square C_{p^2}$. 
\qed

\section{Ordinary Grids}
\label{sec:ordinary}

Since $P_n \square P_n$ has $2n(n-1)$ edges
while $K_{n^2}$ has $n^2 (n^2-1)/2$ edges,
it is clear that if $K_{n^2}$ can be decomposed into
disjoint copies of $P_n \square P_n$,
the number of copies must be $n(n+1)/4$,
so $n$ must be 0 or 3 mod 4.

First we consider the 3-by-3 grid.

\bigskip

{\bf Theorem 3:} 
$K_{9}$ cannot be partitioned into 
three copies of $P_3 \square P_3$.

{\bf Proof:} 
Suppose to the contrary that there were a partition of $K_9$
into subgraphs $G_1$, $G_2$, $G_3$ where each $G_i$ is 
isomorphic to $P_3 \square P_3$ via some isomorphism $\phi_i$.
As before we write the vertices of $K_9$ as 1,\dots,9
and the vertices of $P_3 \square P_3$ as $(i,j)$ with $1 \leq i,j \leq 3$.
Let $v_i$ be the vertex of $K_9$ satisfying $\phi_i (v_i) = (2,2)$
(the middle vertex of $P_3 \times P_3$).
That is, if we adopt the Connections puzzle point of view,
$v_1$, $v_2$, and $v_3$ are the middle elements of the three grids.
For $1 \leq v \leq 9$, let $d_i (v)$ denote the number of elements
of $\{1,\dots,9\}$ that appear next to $v$ in the $i$th grid
(that is, $d_i (v)$ is the degree of $v$ in $G_i$).
Since $d_1 (v_1) + d_2 (v_1) + d_3 (v_1) = 8$
(the degree of $v_1$ in $K_9$)
and since $d_1 (v_1) = 4$
and $d_2 (v_1) \geq 2$ and $d_3 (v_1) \geq 2$,
we must have $d_2 (v_1) = d_3 (v_1) = 2$,
so that $v_1$ appears in a corner cell in both the second and third grids.
Likewise $v_2$ appears in a corner cell in the first and third grids
and $v_3$ appears in a corner cell in the first and second grids.
But this implies that the three elements $v_1, v_2, v_3$
are not adjacent in any of the three grids; this is a contradiction.

\bigskip

We now settle the motivating 4-by-4 example.

\bigskip

{\bf Theorem 4:} 
$K_{16}$ can be partitioned into five copies of $P_4 \square P_4$.

{\bf Proof:} Associate the vertices of $K_{16}$ with 
the elements of $F=\F_2[x]/(x^4+x+1)$ which is a field of order 16.
Note that $x$ is a generator of the multiplicative group $F^*$, 
with $x^{15} = 1$. 
For any choice of $\alpha, \beta, \gamma, \delta\in F$ 
that form a basis for $F$ over $\F_2$,
we can associate the elements of $F$ 
with vertices of $P_4 \square P_4$ as shown in Figure~\ref{fig:grid}.
If the assignment seems random, note that
horizontally-adjacent entries in the middle two columns differ by $\alpha$,
other horizontally-adjacent entries differ by $\beta$,
vertically-adjacent entries in the middle two columns differ by $\delta$,
and other vertically-adjacent entries differ by $\gamma$.
\begin{figure}
\[
\begin{array}{|c|c|c|c|}
\hline
\T\B \beta+\delta & \T\B \delta & \T\B \alpha+\delta & \T\B \alpha+\beta+\delta \\
\hline
\T\B \beta+\gamma+\delta & \T\B \gamma+\delta & \T\B \alpha+\gamma+\delta & \T\B \alpha+\beta+\gamma+\delta \\
\hline
\T\B \beta+\gamma & \T\B \gamma & \T\B \alpha+\gamma & \T\B \alpha+\beta+\gamma \\
\hline
\T\B \beta & \T\B 0 & \T\B \alpha & \T\B \alpha+\beta \\
\hline
\end{array}
\]
\caption{Mapping the 16-element field to the 4-by-4 grid.}
\label{fig:grid}
\end{figure}
Correspondingly let $G(\alpha, \beta, \gamma, \delta)$ 
be the subgraph of $K_{16}$ in which two vertices $v_1, v_2$ 
are joined by an edge if and only if one of the following holds: 
\begin{align*}
    v_1-v_2&=\alpha\quad \mathrm{and}\quad  v_1, v_2\in \mathrm{span}(\alpha, \gamma, \delta)\\
    v_1-v_2&=\beta\\
    v_1-v_2&= \gamma  \\
    v_1-v_2&=\delta\quad \mathrm{and}\quad  v_1, v_2\notin \mathrm{span}(\alpha, \beta, \delta)
\end{align*}
We claim that the graphs $G_i=G(x^{3i}, x^{13+3i}, x^{11+3i}, x^{3+3i})$ 
for $i=0, 1, 2, 3, 4$ form a decomposition of $K_{16}$ 
into 5 copies of $P_4\square P_4$. 
Note that each power of $x^3$ occurs exactly twice 
and everything else occurs just once. 
We get a repeated difference of $x^{3+3i}=v_1-v_2=x^{3(i+1)}$ 
if $v_1, v_2\notin \mathrm{span}(x^{3i}, x^{13+3i}, x^{3+3i})$ and 
$v_1, v_2\in \mathrm{span}(x^{3(i+1)}, x^{11+3(i+1)}, x^{3+3(i+1)}) 
=\mathrm{span}(x^{3+3i}, x^{14+3i}, x^{6+3i})$. 
However,
\begin{align*}
    \mathrm{span}(x^{3i}, x^{13+3i}, x^{3+3i}) 
    &= x^{3i}\ \mathrm{span}(1, x^{13}, x^3)\\
    &=x^{3i}\ \mathrm{span}(1, x^3+x^2+1, x^3)\\
    &=x^{3i}\ \mathrm{span}(x^3, x^3+1, x^3+x^2)\\
    &=x^{3i}\ \mathrm{span}(x^3, x^{14}, x^6)\\
    &=\mathrm{span}(x^{3+3i}, x^{14+3i}, x^{6+3i}).
\end{align*}
Thus there is no repeated difference.
\qed

\bigskip

The five resulting grids are shown in Figure~\ref{fig:five},
with 1 represented by $x^{15}$.
Note that multiplying all the entries in any of the grids 
by $x^3$ in $F$ results in the (cyclically) next of the five grids
(reading left to right in the first row, then down,
then right to left in the second row, then up).
\begin{figure}
\begin{gather*}
\begin{array}{|c|c|c|c|}
\hline
\T\B x^8 & \T\B x^3 & \T\B x^{14} & \T\B x^2 \\
\hline
\T\B x^7 & \T\B x^5 & \T\B x^{10} & \T\B x^9 \\
\hline
\T\B x^4 & \T\B x^{11} & \T\B x^{12} & \T\B x \\
\hline
\T\B x^{13} & \T\B 0 & \T\B x^{15} & \T\B x^6 \\
\hline
\end{array}
\ \ \ \ 
\begin{array}{|c|c|c|c|}
\hline
\T\B x^{11} & \T\B x^6 & \T\B x^2 & \T\B x^5 \\
\hline
\T\B x^{10} & \T\B x^8 & \T\B x^{13} & \T\B x^{12} \\
\hline
\T\B x^7 & \T\B x^{14} & \T\B x^{15} & \T\B x^{4} \\
\hline
\T\B x^{1} & \T\B 0 & \T\B x^3 & \T\B x^9 \\
\hline
\end{array}
\ \ \ \ 
\begin{array}{|c|c|c|c|}
\hline
\T\B x^{14} & \T\B x^9 & \T\B x^5 & \T\B x^8 \\
\hline
\T\B x^{13} & \T\B x^{11} & \T\B x^{1} & \T\B x^{15} \\
\hline
\T\B x^{10} & \T\B x^{2} & \T\B x^{3} & \T\B x^{7} \\
\hline
\T\B x^{4} & \T\B 0 & \T\B x^6 & \T\B x^{12} \\
\hline
\end{array} 
\\\\
\begin{array}{|c|c|c|c|}
\hline
\T\B x^{2} & \T\B x^{12} & \T\B x^8 & \T\B x^{11} \\
\hline
\T\B x^{1} & \T\B x^{14} & \T\B x^{4} & \T\B x^{3} \\
\hline
\T\B x^{13} & \T\B x^{5} & \T\B x^{6} & \T\B x^{10} \\
\hline
\T\B x^{7} & \T\B 0 & \T\B x^9 & \T\B x^{15} \\
\hline
\end{array}
\ \ \ \ 
\begin{array}{|c|c|c|c|}
\hline
\T\B x^{5} & \T\B x^{15} & \T\B x^{11} & \T\B x^{14} \\
\hline
\T\B x^{4} & \T\B x^{2} & \T\B x^{7} & \T\B x^{6} \\
\hline
\T\B x^{1} & \T\B x^{8} & \T\B x^{9} & \T\B x^{13} \\
\hline
\T\B x^{10} & \T\B 0 & \T\B x^{12} & \T\B x^{3} \\
\hline
\end{array}
\end{gather*}
\caption{A $P_4 \square P_4$ grid design.}
\label{fig:five}
\end{figure}
The reader can check that every pair of elements $u,v \in F$ 
occurs adjacently in exactly one of the five grids.

\section{Comments} \label{sec:comments}

The simplest cases we were unable to resolve are
$C_{15} \square C_{15}$ and $P_7 \square P_7$.
We think that the case of $P_p \square P_p$
with $p$ congruent to 0 mod 4
might yield to some variant of the construction
presented in section~\ref{sec:ordinary}.

One might look at higher-power grids
such as $C_n \square C_n \square \cdots C_n$
and $P_n \square P_n \square \cdots P_n$.
It would also be natural to consider products
like $C_m \square C_n$ or $P_m \square P_n$
or for that matter $C_m \square P_n$.

\bigskip

\noindent
{\sc Acknowledgment:}
We thank Erik Demaine for hacking into the New York Times website's software
to determine that it does indeed use a pseudorandom number generator.

\end{document}